\newcommand{\ZZ}{\mathbb{Z}}
\title{Universal Cycles on 3--Multisets}
\author{
Joshua $\ZZ$ahl
\thanks{ Caltech, Pasadena, CA 91125, USA. Email: jzahl@zahl.ca.
Research supported in part by NSF grant 0552730.} \and Tobias L.
Johnson
\thanks {Yale, New Haven, CT 06520, USA. Email: tobias.l.johnson@gmail.com. Research supported in part by NSF grant
0552730.} }

\date{}

\documentclass [11pt]{article}
\usepackage{amsfonts}
\usepackage{amsmath}
\usepackage{amsthm}

\usepackage{units}

\begin{document}
\maketitle



\begin{abstract}
\noindent Consider the collection of all $t$--multisets of
$\{1,\ldots, n\}$. A \emph{universal cycle on multisets} is a
string of numbers, each of which is between $1$ and $n$, such that
if these numbers are considered in $t$--sized windows, every
multiset in the collection is present in the string precisely
once. The problem of finding necessary and sufficient conditions
on $n$ and $t$ for the existence of universal cycles and similar combinatorial structures was first addressed by
DeBruijn in 1946 (who considered $t$--tuples instead of
$t$--multisets). The past 15 years has seen a resurgence of
interest in this area, primarily due to Chung, Diaconis, and
Graham's 1992 paper on the subject. For the case $t=3$, we
determine necessary and sufficient conditions on $n$ for the
existence of universal cycles, and we examine how this technique
can be generalized to other values of $t$.
\end{abstract}
\maketitle
\section{Introduction and Previous Work}
Consider the collection of all $t$--multisets over the universe
$[n]=\{1,\ldots, n\}$. A universal cycle (ucycle) on multisets is
a cyclic string $X=a_1a_2...a_k$ with $a_i\in[n]$ for which the
collection $\big\{\{a_1,a_2,...,a_t\},$
$\{a_2,a_3,...,a_{t+1}\},...,$ $\{a_{k-t+1},a_{k-t+2},...,a_k\},$
$\{a_{k-t+2},a_{k-t+3},...,a_k,a_1\},...
\{a_k,a_1,a_2,...,a_{t-1}\} \big\}$ is precisely the collection of
all $t$--multisets over $[n]$, i.e. each $t$--multiset over $[n]$
occurs precisely once in the above collection. For the remainder
of this paper, the term universal cycle will refer to universal
cycles on multisets unless noted otherwise.
%

Universal cycles do not exist for every value of $n$ and $t$.
Indeed, simple symmetry arguments show that each of the numbers
$1,\ldots,n$ must occur an equal number of times in the ucycle.
Since the length of the ucycle is equal to the number of
$t$--multisets over $[n]$, which is $\binom{n+t-1}{t}$, we must
have that $n\!\left.|\binom{n+t-1}{t}\right..$ While this
condition is necessary, it is not sufficient for the existence of
ucycles.

To date, the bulk of research on ucycles has been devoted to
studying ucycles over sets (as opposed to multisets). Ucycles over
sets are constructed in the same fashion as ucycles over
multisets, except that we consider the collection of all $t$--sets
over $[n]$ instead of the collection of all $t$--multisets, and
our divisibility condition becomes $n\!\left|\binom{n}{t}\right.$.
In \cite{Chung}, Chung, Diaconis and Graham conjectured that for
each value of $t$, there exists a number $n_0(t)$ such that
universal cycles exist for $n\left|\binom{n}{t}\right.$ and $n\geq
n_0(t)$. In \cite{Hurlbert}, Hurlbert consolidated and extended
previous work, verifying the conjecture for $t=2$ and $3$ and
developing partial results for $t=4$ and $6$. In \cite{Godbole},
Godbole et al.~considered universal cycles over multisets for the
case $t=2$, and verified the analogous form of the Chung--Diaconis--Graham conjecture (i.e. with the modified divisibility
criterion) for this case. This work is of particular interest
because Godbole et al.~used a new inductive technique to arrive at
their proof, and in this paper we extend this technique to the
case $t=3$. This new work suggests that the inductive method is a
promising way of addressing the Chung--Diaconis--Graham
conjecture. We also consider a second proof of the conjecture for
the $t=3$ multiset case, which builds off universal cycles on sets
and lends itself more easily to generalization.
\section{An Inductive Proof for Universal Cycles on 3--Multisets}
For $t=3$, the condition $n\left|\binom{n+t-1}{t}\right.$ implies
that $t\equiv 1$ or $2$ (mod 3). We will consider the case
$n\equiv 1$ (mod 3), as the other case can be dealt with
similarly. We will show that for $n\geq 4$, universal cycles exist
whenever $n$ satisfies $n\left|\binom{n+t-1}{t}\right.$

Before describing the proof itself, we will define some
terminology that will be useful for describing universal cycles.
We say that a cyclic string $X=a_1a_2...a_k$ \emph{contains} the
multiset collection $\mathcal{I}$ if
$\mathcal{I}=\big\{\{a_1,a_2,a_3\},\{a_2,a_3,a_4\},...,\{a_{k-2},a_{k-1},{a_k}\},\{a_{k-1},a_k,a_1\},\{a_k,a_1,a_2\}
\big\}$, where each of these sets must be distinct. Clearly
$k=\binom{n+2}{3}$, since this is the number of $3$--multisets on
$[n]$.

For a string $X=a_1a_2...a_k$, we call the \emph{lead-in} of $X$
the substring $a_1a_2$ and the \emph{lead-out} the substring
$a_{k-1}a_k$.

%
Now, consider the collection of all 3--multisets over $[n]$. We
shall partition this collection into four subcollections. Let
$\mathcal{A}$ be the collection of all 3--multisets over $[n-3]$,
and let $\mathcal{B}$ be the collection of all $3$--multisets over
$\{n-2,n-1,n\}$ and $[n-6]$ which contain at least one element
from $\{n-2,n-1,n\}$. Let $\mathcal{C}$ be the collection of all
3--multisets with one or two elements from $\{n-5,n-4,n-3\}$ and
one or two elements from $\{n-2,n-1,n\}$, and let $\mathcal{D}$ be
the collection of all 3--multisets with one element from each of
$[n-6]$, $\{n-5,n-4,n-3\}$, and $\{n-2,n-1,n\}$. We can see that
$\mathcal{A,B,C},$ and $\mathcal{D}$ are disjoint, and that their
union is the collection of all 3--multisets on $[n]$, as desired.

Now, let $S$ be a universal cycle on $[n-6]$, and since $1,1,1$
must occur somewhere in $S$ and the beginning of $S$ is arbitrary,
we shall have $S$ begin with $1,1,1$. We shall also select $S$ so
that its lead-out is $n-6,n-7$. Thus $S$, when considered as a
cyclic string, contains all $3$--multisets over $[n-6]$, and when
considered as a non-cyclic string, contains all $3$--multisets
except $\{1,n-7,n-6\}$ and $\{1,1,n-7\}$. Let $T$ be a string over
$[n-3]$ such that $ST$---the concatenation of $S$ and $T$---is a
universal cycle over $[n-3]$. It is not clear that such a $T$ must
exist, but we shall find a specific example shortly. In the
example we will find, $T$ will begin with $1,1$ and will end with
$n-3,n-4$. Since $T$ begins with $1,1$, the string $ST$ contains
the multisets $\{1,n-7,n-6\},\ \{1,1,n-7\}$. We can see that the
cyclic string $ST$ contains all of the multisets in $\mathcal{A}$,
and that when $ST$ is considered as a non-cyclic string, it
contains $\mathcal{A}\backslash\big\{\{1,n-4,n-3\},\
\{1,1,n-4\}\big\}$. Now, consider the string $T^\prime$ obtained
by taking $T$ and replacing each instance of $n-5$ by $n-2$, $n-4$
by $n-1$, and $n-3$ by $n$. Since $T$ contained all multisets over
$[n-3]$ which contained at least one element from
$\{n-5,n-4,n-3\}$, we have that $T^\prime$ contains all multisets
over $\{n-2,n-1,n\}$ and $[n-6]$ which contain at least one
element from $\{n-2,n-1,n\}$, i.e. $T^\prime$ contains all the
multisets in $\mathcal{B}$. Since the lead-in of $T$ is $1,1$, the
lead-in of $T^\prime$ is also $1,1$, and since $T$ ends with
$n-3,n-4$, $T^\prime$ ends with $n,n-1$. If we consider the cyclic
string $STT^\prime$, we can see that this string contains all the
multisets in $\mathcal{A}\cup \mathcal{B}$, while the non-cyclic
version of this string is missing the multisets $\{1,n-1,n\},\
\{1,1,n-1\}$.

For notational convenience, we will use the following assignments:
$a:=n-5,\ b:=n-4,\ c:=n-3,\ d:=n-2,\ e:=n-1,$ and $f:=n$. Now,
consider the following string:

\begin{eqnarray*}
V=&&\mathrm{be}(n-6)\mathrm{af}(n-7)\mathrm{be}(n-8)\mathrm{af}(n-9)...\mathrm{af}1\mathrm{be}\\
&&\ \ \ \mathrm{ad}(n-6)\mathrm{ce}(n-7)\mathrm{ad}(n-8)\mathrm{ce}(n-9)...\mathrm{ce}1\mathrm{ad}\\\
&&\ \ \ \ \ \
\mathrm{cf}(n-6)\mathrm{bd}(n-7)\mathrm{cf}(n-8)\mathrm{bd}(n-9)...\mathrm{bd}1\mathrm{cfe}.
\end{eqnarray*}
We can see that this string contains every multiset in
$\mathcal{D}$, as well as the multisets $\{a,b,e\}$, $\{a,d,e\}$,
$\{a,c,d\}$, and $\{c,d,f\}$. Now, the following string (found
with the aid of a computer) contains all of the multisets in
$\mathcal{C}\backslash \big\{\{a,b,e\}$, $\{a,d,e\}$, $\{a,c,d\}$,
$\{c,d,f\} \big\}$:
$$U=\mathrm{aaffc\phantom{1}aeebb\phantom{1}decec\phantom{1}bddcc\phantom{1}fbada\phantom{1}dfbf}.$$

Note that while the multisets $\{b,b,f\}$ and $\{b,e,f\}$ are not
present in the above string $U$, they are present in the
concatenation of $U$ with $V$. Similarly, while $U$ does not
contain $\{a,e,f\}$ and $\{a,a,f\}$, these multisets are present
in the concatenation of $T^\prime$ with $U$.

Now, we can see that the string $STT^\prime UV$ is a universal
cycle over $[n]$ because the non-cyclic string $STT^\prime$
contained all the multisets in
$\mathcal{A}\cup\mathcal{B}\backslash\big\{\{1,n-1,n\},\
\{1,1,n-1\} \big\}$, and it is precisely the multisets
$\{1,n-1,n\}$ and $\{1,1,n-1\}$ which are obtained by the
wrap-around of the lead-out of $V$ with the lead-in of $S$. The
lead-in and lead-out of the other strings has been engineered so
as to ensure that each multiset occurs precisely once.

This completes the induction proof, since the string $ST$ is a
universal cycle over $[n-3]$ (taking the place of $S$ in the
previous iteration of the induction), and the string $T^\prime UV$
extends this cycle to $[n]$ (taking the place of $T$ in the
previous iteration of the induction). Also note that $T^\prime UV$
begins with $1,1$ and ends with $n,n-1$, as required for the
induction hypothesis.

Thus, all that remains is the find a base case from which the
induction can proceed. A possible base case (there are many) for
$n-6=4,\ n-3=7$ is
\begin{eqnarray*}S&=&11144\phantom{1}42223\phantom{1}33121\phantom{1}24343\\
T&=&11522\phantom{1}63374\phantom{1}45166\phantom{1}27732\phantom{1}57366\phantom{1}77135\phantom{1}34641\phantom{1}71555\phantom{1}36127\phantom{1}42556\phantom{1}66477\phantom{1}75526\phantom{1}4576,\end{eqnarray*}
which would lead to
\begin{eqnarray*}
T^\prime&=&11822\phantom{1}93304\phantom{1}48199\phantom{1}20032\phantom{1}80399\phantom{1}00138\phantom{1}34941\phantom{1}01888\phantom{1}39120\phantom{1}42889\phantom{1}99400\phantom{1}08829\phantom{1}4809\\
U&=&55007\phantom{1}59966\phantom{1}89797\phantom{1}68877\phantom{1}06585\phantom{1}8060\\
V&=&69450\phantom{1}36925\phantom{1}01695\phantom{1}84793\phantom{1}58279\phantom{1}15870\phantom{1}46837\phantom{1}02681\phantom{1}709,
\end{eqnarray*}
Where ``0'' denotes 10 and the spacings have been added to
increase readability.

All of the work up to this point has dealt with $n \equiv 1 $ (mod
3). The proof for $n\equiv2$ (mod 3) is similar, so it has been
omitted for the sake of brevity.
\section{A Second Proof of the Existence of Ucycles on 3--Multisets}
In this proof, we construct a ucycle on 3--multisets of $[n]$ by modifying a ucycle on 3--subsets of $[n]$.  (We know from \cite{Hurlbert} that ucycles on 3--subsets of $[n]$ exist for all $n \geq 8$ not divisible by $3$.)  Before giving the proof, we introduce two terms.  We call each element of $[n]$ a \emph{letter}, and each $a_i$ in the ucycle $X=a_1\ldots a_k$ a \emph{character}.  To summarize, a ucycle on 3--multisets of $[n]$ is made up of $\binom{n+t-1}{t}$ characters, each of which equals one of $n$ letters.

To demonstrate the proof's technique, we will first use an argument similar to it to create ucycles on 2--multisets from ucycles on 2--subsets.
We start with this ucycle on 2--subsets of $[5]$:
$$
1234513524
$$
Then, we repeat the first instance of every letter to create the following ucycle on 2--multisets:
$$
112233445513524
$$
The technique works because repeating a character $a_i$ as above adds the multiset $\{a_i,a_i\}$ to the ucycle and has no other effect.

To use this technique on ucycles on 3--subsets, we repeat not single characters, but pairs of characters.  For example, changing
$$
\ldots a_{i-1}a_ia_{i+1}a_{i+2}\ldots
$$
to
$$
\ldots a_{i-1}a_ia_{i+1}a_ia_{i+1}a_{i+2} \ldots
$$
has only has the effect of adding the 3--multisets $\{a_i,a_i,a_{i+1}\}$ and $\{a_i,a_{i+1},a_{i+1}\}$ to the cycle.  In order to use this technique, we will need to know which consecutive pairs of letters appear in a ucycle on 3--subsets.  For instance, the following ucycle (generated using methods from \cite{Hurlbert}) on 3--subsets of $[8]$ contains every unordered pair of letters as consecutive characters but $\{1,5\}$, $\{2,6\}$, $\{3,7\}$, and $\{4,8\}$:
$$
1235783\ 6782458\ 3457125\ 8124672\ 5671347\ 2346814\ 7813561\ 4568236
$$
(The spaces in the cycle are added only for readability.)
This ucycle is missing 4 pairs, which happens to be $n/2$.  This is no coincidence: in fact, this is the most pairs that a ucycle on 3--subsets can fail to contain.

\newtheorem*{missingpairs}{Lemma}
\begin{missingpairs}
No two unordered pairs not appearing as consecutive characters in a ucycle on 3--subsets have a letter in common.  A ucycle can hence be missing at most $n/2$ pairs of letters.
\end{missingpairs}
\begin{proof}
Suppose that we have a ucycle on 3--subsets that contains neither $a$ and $b$ as consecutive characters, nor $a$ and $c$ as consecutive characters, where $a,b,c \in [n]$.  Then the ucycle does not contain the 3--subset $abc$, for all permutations of $abc$ contain either $a$ and $b$ consecutively, or $a$ and $c$ consecutively.  But this is a contradiction, as a ucycle by definition contains all 3--subsets.

Hence, no two pairs of characters missing in the ucycle can have a letter in common.  By the pigeonhole principle, the ucycle can be missing at most $n/2$ pairs of letters.
\end{proof}

With this lemma, we can finish our proof, creating a ucycle on 3--multisets of $[n]$ whenever $n$ is not divisible by 3.
First, we consider the case when $n$ is even.  Let $X$ be a ucycle on 3--subsets of $[n]$.
Let $x_1,\ldots,x_n$ be a permutation of $[n]$ such that
\begin{itemize}
\item $x_1$ equals the first character in $X$.
\item $x_n$ equals the last character in $X$.
\item The list $\{x_1,x_2\}, \{x_3,x_4\}, \ldots, \{x_{n-1}x_n\}$ contains all unordered pairs of letters not contained as consecutive characters in $X$, which is possible by our lemma.  (If $X$ is missing exactly $n/2$ pairs of letters, these pairs will be exactly the pairs missing from $X$.  If $X$ is missing fewer than $n/2$ pairs of letters, then the pairs consist of all missing pairs of letters, plus the remaining letters paired arbitrarily.)
\end{itemize}
Make $X'$ by repeating the first instance of every
unordered pair of letters in $X$ except for $\{x_1,x_2\},
\{x_2,x_3\},$ $\ldots,$ $\{x_{n-1},x_n\},\{x_n,x_1\}$.  The cycle $X'$
now contains all multisets except 
\setlength\arraycolsep{1pt}
$$
\{x_1,x_1,x_1\},\ldots,\{x_n,x_n,x_n\}\\
$$
$$
\{x_1,x_1,x_2\},\{x_1,x_2,x_2\},\{x_2,x_2,x_3\},\{x_2,x_3,x_3\},\ldots,\{x_n,x_n,x_1\},\{x_n,x_1,x_1\}
$$
Now, add the string $x_1x_1x_1x_2x_2x_2\ldots x_nx_nx_n$ to the
end of $X'$ to create $X''$.  This provides exactly the missing multisets, creating a
ucycle on 3--multisets.

For example, when $n=8$, we start with the following ucycle on 3--subsets:
\begin{eqnarray*}
X & = & 1235783\ 6782458\ 3457125\ 8124672\\
&& 5671347\ 2346814\ 7813561\ 4568236\\
\end{eqnarray*}
The ucycle on 3--subsets $X$ does not contain the pairs $\{1,5\}$, $\{2,6\}$, $\{3,7\}$, and $\{4,8\}$.  Hence, we set
\begin{eqnarray*}
x_1 & = & 1,\ x_2=5,\ x_3=3,\ x_4=7\\
x_5 & = & 4,\ x_6=8,\ x_7=2,\ x_8=6
\end{eqnarray*}
Note that $x_1$ equals the first character of $X$, and $x_8$ equals the last.

Now, we repeat the first instance of every unordered pair except for $\{1,5\}$, $\{5,3\}$, $\{3,7\}$, $\{7,4\}$, $\{4,8\}$, $\{8,2\}$, $\{2,6\}$, and $\{6,1\}$.  (Note that four of these pairs do not appear in $X$.  If some of these pairs actually did appear in $X$, because $X$ was missing fewer than $n/2$ pairs of letters, it would not affect the proof.):
\begin{eqnarray*}
X' & = &12123235757878383\
63676782424545858\
3434571712525\
81812464672\\
&& 56567131347\
2723468681414\
7813561\
4568236
\end{eqnarray*}
Finally, we add the string $x_1x_1x_1\ldots x_nx_nx_n$ to complete the ucycle:
\begin{eqnarray*}
X'' & = &12123235757878383\
63676782424545858\
3434571712525\
81812464672\\
&& 56567131347\
2723468681414\
7813561\
4568236\\
&& 111555333777444888222666
\end{eqnarray*}

The proof is similar when $n$ is odd, and we omit it for
the sake of brevity.
\section{Further Directions and Remarks}
Both of the proofs given above suggest natural extensions to the
$t=4$ and larger cases, and it is simple to use the
techniques described above to create a proof sketch. In personal
correspondence, Glenn Hurlbert indicated that his technique for creating ucycles on sets in \cite{Hurlbert} can also be used to create ucycles on multisets.
Though this provides a more concise proof for the existence of ucycles on 3--multisets, the two proofs presented may prove useful by their introduction of new techniques for approaching ucycles.  The first proof is notable for its use of induction, a technique which has not been used before to create ucycles.  The second, while it is tied to ucycles on sets, is not tied to any particular approach for creating ucycles on sets; it could perhaps be extended to situations to which Hurlbert's technique cannot.

For values of $n$ and $t$ for which ucycles do exist, one
interesting question is how many ucycles exist. Clearly each
ucycles has $n!$ representations, since there are $n!$
permutations of $1,\ldots,n$. However, when searching for ucycles
using a computer, vast numbers of \emph{distinct} (i.e. not
differing merely by a permutation of $1,\ldots,n$) ucycles were
found. Currently, it is not clear whether $N(n,t)$, the number of
distinct ucycles for a given value of $n$ and $t$, is a function
that has a simple description.
\bibliographystyle{amsplain}

\end{document}